\documentclass[reqno,11pt]{amsart} 
\usepackage{amsmath}
\usepackage{amssymb}
\usepackage{amsthm}
\usepackage{verbatim}
\usepackage{xcolor}
\usepackage{graphics}

\newcommand\NoBlackBoxes{\global\overfullrule0pt}

\NoBlackBoxes
\theoremstyle{plain}

\begin{document}

\title{ON SUBGRADIENTS OF CONVEX FUNCTIONS AND ORLICZ \\
PSEUDO-NORMS FOR VECTOR-VALUED FUNCTIONS
}

\author{Sergey G. Bobkov$^{1}$}
\thanks{1) 
School of Mathematics, University of Minnesota, Minneapolis, MN, USA,
bobkov@math.umn.edu. 
}

\author{Friedrich G\"otze$^{2}$}
\thanks{2) Faculty of Mathematics,
Bielefeld University, Germany,
goetze@math-uni.bielefeld.de.
}
\thanks{3) Research supported by the NSF grant DMS-2154001 and
	the GRF – SFB 1283/2 2021–317210226}
	

\begin{abstract}
We discuss variants of construction of measurable subgradients for multivariate convex functions 
and the  problem of characterization of the $\Delta_2$-condition in terms of their directional
derivatives. Furthermore we study related basic properties of Luxemburg and Orlicz pseudo-norms for 
vector-valued functions.
\end{abstract}

\maketitle
\markboth{Sergey G. Bobkov and Friedrich G\"otze}{Orlicz pseudo-norms}

\def\theequation{\thesection.\arabic{equation}}
\def\E{{\mathbb E}}
\def\R{{\mathbb R}}
\def\C{{\mathbb C}}
\def\P{{\mathbb P}}
\def\Z{{\mathbb Z}}
\def\L{{\mathbb L}}
\def\T{{\mathbb T}}

\def\G{\Gamma}

\def\Ent{{\rm Ent}}
\def\var{{\rm Var}}
\def\Var{{\rm Var}}

\def\H{{\rm H}}
\def\Im{{\rm Im}}
\def\Tr{{\rm Tr}}
\def\s{{\mathfrak s}}

\def\k{{\kappa}}
\def\M{{\cal M}}
\def\Var{{\rm Var}}
\def\Ent{{\rm Ent}}
\def\O{{\rm Osc}_\mu}

\def\ep{\varepsilon}
\def\phi{\varphi}
\def\vp{\varphi}
\def\F{{\cal F}}

\def\be{\begin{equation}}
\def\en{\end{equation}}
\def\bee{\begin{eqnarray*}}
\def\ene{\end{eqnarray*}}

\thispagestyle{empty}

\vskip5mm
\section{{\bf Introduction}}
\setcounter{equation}{0}

\vskip2mm
This note is devoted to the study of several general questions about convex cost functions 
on $\R^n$ related to the construction of subgradients, characterizations of the $\Delta_2$-condition, 
and some treatment of Orlicz seminorms in the multivariate setting.
We found that these questions for general convex functions were not properly addressed in 
the literature but turned out to be  crucial for  our current work  on energy-type estimates 
in transport problems.  In this upcoming work we intend to extend  bounds of M. Ledoux \cite{L}  
for the Kantorovich transport distances $W_p$ to transport distances based on general convex cost 
functions. We think that this note might be 
of independent interest for further transport problems in $\R^n$ as well. 

The first part of this paper is devoted to the construction of subgradients of convex
functions on $\R^n$ which is a standard topic in the area of non-smooth optimization. 
As it turns out, the existence of measurable subgradients is required
in order to establish in a rigorous way the relationship 
between the Luxemburg and the Orlicz pseudo-norms of functions with values in $\R^n$.
This relationship, which we consider in the third part of the paper,
extends a known result in the scalar case to the multidimensional setting. 
In the second part, we discuss the $\Delta_2$-condition.

\vskip7mm
\section{{\bf Differentiability and Subgradients of Convex Functions}}
\setcounter{equation}{0}

\vskip2mm
\noindent
To start with, let us recall several definitions and standard results
about differentiability and subgradients of convex functions.

Let $L:\R^n \rightarrow \R$ be a convex function.
It is well-known that $L$ is locally Lipschitz, hence almost everywhere differentiable, by
the Rademacher theorem (cf. \cite{H} for refinements). Thus, the  set $E$ of all points 
$x \in \R^n$ where $L$ is differentiable has a full Lebesgue measure.
It is also known that $x \in E$ if and only if the partial derivatives
$$
\partial_i L(x) = \lim_{\ep \rightarrow 0} \frac{L(x + \ep e_i) - L(x)}{\ep}, \quad i = 1,\dots,n,
$$
exist, where $\{e_i\}_{i=1}^n$ denotes the canonical basis in $\R^n$.

Anyway, since the functions $(L(x + \ep e_i) - L(x))/\ep$ are monotone in $\ep$,
the right and left partial derivatives
\bee
\partial_i^+ L(x) 
 & = &
\lim_{\ep \downarrow 0} \frac{L(x + \ep e_i) - L(x)}{\ep}, \\
\partial_i^- L(x) 
 & = &
\lim_{\ep \uparrow 0} \frac{L(x + \ep e_i) - L(x)}{\ep}
\ene
exist, are finite, and represent Borel measurable functions satisfying
$\partial_i^- L(x) \leq \partial_i^+ L(x)$. As a consequence, the set $E$ may be represented as
$$
E = \{x \in \R^n: \partial_i^+ L(x) = \partial_i^- L(x)\, \ {\rm for \ any} \ i = 1,\dots,n\}.
$$
Hence this set is Borel measurable, and the gradient function $\nabla L(x)$ is Borel measurable as well
on $E$. In fact, the function $\nabla L$ is continuous on $E$ (cf. \cite{R}). This implies
that if $L$ is differentiable everywhere, then it is $C^1$-smooth. For various results about
convex functions, let us also refer to the book \cite{Ho}.

For $x \in \R^n$, consider the subdifferential
\be
\partial L(x) = \Big\{y \in \R^n: L(z) \geq L(x) + \left<y,z-x\right> \ {\rm for \ all} \ z \in \R^n\Big\},
\en
which is not empty, by the convexity of $L$. Vectors in this set are called subgradients of $L$
at the point $x$. If $L$ is differentiable at $x$, then  $\partial L(x)$ contains only one vector $y = \nabla L(x)$.

In the general situation, $\partial L(x)$ is convex and compact. Convexity and closeness
is obvious. Since for $y \in \partial L(x)$, 
$$
\left<y,h\right> \leq L(x+h) - L(x) \ \ {\rm for \ all} \ h \in \R^n,
$$
we have
\be
|y| \, = \, \frac{1}{r}\,\sup_{|h| \leq r} \left<y,h\right> \, \leq \, 
\frac{1}{r}\,\sup_{|h| \leq r} (L(x+h) - L(x)), \quad r>0.
\en
So, the set $\partial L(x)$ is bounded, with diameter depending on $x$.

As any convex compact set, the subdifferential at a fixed point $x \in \R^n$ 
may be characterized by the support function
$$
h_{\partial L(x)}(\theta) = \sup_{y \in \partial L(x)} \left<y,\theta\right>, \quad \theta \in \R^n.
$$
For this aim, one should involve the directional derivatives
\begin{eqnarray}
L'(x,\theta) 
 & = & 
\lim_{\ep \downarrow 0} \frac{L(x + \ep \theta) - L(x)}{\ep} \nonumber \\
 & = & 
\inf_{\ep > 0} \frac{L(x + \ep \theta) - L(x)}{\ep}, \quad x,\theta \in \R^n.
\end{eqnarray}
Similarly to the case of right derivatives, this limit exists, is finite, and represents a Borel measurable function
in two variables $(x,\theta)$. In fact, $L'$ is upper semi-continuous on $\R^n \times \R^n$
(cf. \cite{R}, Corollary 24.5.1). Note also that the function $\theta \rightarrow L'(x,\theta)$ 
is positive homogeneous and convex. Returning to the support function,
it may be shown (cf. \cite{R}, Theorem 23.4) that, for all $x,\theta \in \R^n$,
$$
h_{\partial L(x)}(\theta) = L'(x,\theta).
$$

Now, consider the conjugate function (also called the dual and Legendre transform).
\be
L^*(x) = \sup_{y \in \R^n} \big[\left<x,y\right> - L(y)\big], \quad x \in \R^n.
\en
It is a convex function with values in $(-\infty,\infty]$, and $L^{**} = L$.
If $L$ satisfies the super-linear growth condition
$$
L(x)/|x| \rightarrow \infty \ \ {\rm as}\, \ |x|\rightarrow \infty,
$$
then the conjugate functions is finite for all $x$ and possesses a similar growth property.

The dual transform may be used in the following well-known characterization
(cf. \cite{R}, Theorem 23.5), which we complement by an upper bound for the conjugate
$L^*$ on the subdifferential of $L$.

\vskip5mm
{\bf Proposition 2.1.} {\sl Given $x \in \R^n$, the value $L^*(y)$ is finite whenever $y \in\partial L(x)$.
Moreover,
\be
\left<x,y\right> = L(x) + L^*(y) \ \Longleftrightarrow \ y \in\partial L(x).
\en
As a consequence, for any $r>0$
\be
\sup_{y \in \partial L(x)} L^*(y) \, \leq \, \frac{|x|}{r}\,\sup_{|h| \leq r} \big(L(x+h) - L(x)\big) - L(x).
\en
}

\vskip3mm
{\bf Proof.} According to the definition (2.1), the property $y \in \partial L(x)$ is equivalent to
\be
\left<x,y\right> - L(x) \geq \sup_z \big[\left<z,y\right> - L(z)\big].
\en
But the last supremum is just $L^*(y)$, which has to be bounded from above by the left hand-side.
Hence it is finite. On the other hand, the left-hand side does not exceed $L^*(y)$,
according to the definition (2.4). This shows that $\left<x,y\right> - L(x) = L^*(y)$.
Conversely, the latter equality implies (2.7), by the very definition of the conjugate function.

Finally, combining the equality in (2.5) with the upper bound (2.2), we obtain (2.6):
\bee
\sup_{y \in \partial L(x)} L^*(y) 
 & = &
\sup_{y \in \partial L(x)} \left<x,y\right> - L(x) \\
 & \leq & 
|x|\,\sup_{y \in \partial L(x)} |y| - L(x) \\
 & \leq & 
\frac{|x|}{r}\,\sup_{|h| \leq r} \big(L(x+h) - L(x)\big) - L(x).
\ene
\qed

\vskip5mm
An immediate consequence of (2.5) is the duality
$
y \in\partial L(x) \Longleftrightarrow x \in\partial L^*(y).
$

\vskip7mm
\section{{\bf Existence of Measurable Subgradients}}
\setcounter{equation}{0}

\vskip2mm
\noindent
We will need to construct a subgradient function $y = y(x)$ as a map depending on $x$ in a measurable way.
As we will see, this question appears naturally in the study of the relationship between Orlicz-type norms.

\vskip5mm
{\bf Proposition 3.1.} {\sl  For any convex function $L$ on $\R^n$, there exists a Borel measurable map 
$y:\R^n \rightarrow \R^n$ such that $y(x) \in \partial L(x)$ for any $x \in \R^n$. 
}

\vskip5mm
As a natural candidate for such a map, one may consider the barycenter
of the subdifferential $\partial L(x)$,
\be
y(x) = \int_{\partial L(x)} z\,dm_x(z) = \frac{1}{|\partial L(x)|} \int_{\partial L(x)} z\,dz.
\en
Here $m_x$ denotes a uniform distribution on $\partial L(x)$, that is, the Lebesgue measure on the 
smallest affine subspace $E(x)$ of $\R^n$ containing $\partial L(x)$, restricted to this subdifferential and 
normalized, where $|\partial L(x)|$ is its $k$-dimensional volume, and $k$ is the dimension of $E(x)$.
In particular, if $L$ is differentiable at $x$, then
$\partial L(x) = \{\nabla L(x)\}$ and $m_x = \delta_{\nabla L(x)}$ is the delta-measure at the point $x$.

By construction, $y(x) \in \partial L(x)$ for any $x \in \R^n$. 
It remains to show that this map is Borel measurable. To prove this, consider
the space $\mathbb K_n$ of all non-empty compact subsets $K$ of $\R^n$.
It becomes a complete separable metric space with the Hausdorff metric
$$
\rho(K_1,K_2) = \inf \big\{\ep \geq 0: K_1 \subset K_2 + \ep B, \ K_2 \subset K_1 + \ep B\big\},
$$
where $B$ denotes the unit Euclidean ball in $\R^n$. Here we use the Minkowski summation
$$
a K + b L = \{ax+by: x \in K, \ y \in L\}
$$ 
for non-empty subsets $K,L$ of $\R^n$ and scalars $a,b$.

Let $(\Omega,d)$ be a metric space. 
Given a map $T:\Omega \rightarrow {\mathbb K}_n$, its continuity is understood with 
respect to $\rho$. We say that $T$ is upper semi-continuous, if
\be
\forall\,\ep>0 \ \, \forall x_0 \in \Omega \ \, \exists\, \delta>0 \ 
\Big[ x \in \Omega, \, d(x,x_0) < \delta \, \Longrightarrow \, T(x) \subset T(x_0) + \ep B \Big].
\en
This is equivalent to the set-theoretical upper semi-continuity
$$
\limsup_{x \rightarrow x_0}\, T(x) \, \equiv \, 
\bigcap_{\delta > 0}\, \bigcup_{d(x,x_0) < \delta} T(x) \, = \, T(x_0).
$$
Indeed, from (3.2), it follows that $\limsup_{x \rightarrow x_0} T(x) \subset T(x_0) + \ep B$ 
for any fixed $\ep>0$. But, 
$$
\bigcap_{\ep>0}\, \big(T(x_0) + \ep B\big) = T(x_0), 
$$
by the compactness of $T(x_0)$.
We will need the following.

\vskip5mm
{\bf Lemma 3.2.} {\sl Any upper semi-continuous map $T$ on the metric space $(\Omega,d)$
with values in ${\mathbb K}_n$ is Borel measurable.
}

\vskip5mm
{\bf Proof.} The elements $K$ of $\mathbb K_n$ may be identified
with their support functions
$$
h_K(x) = \sup_{y \in K} \left<x,y\right>, \quad x \in \R^n.
$$
Any support function is positive homogeneous and convex on $\R^n$. This operation preserves 
inclusion: $K_1 \subset K_2$ if and only if $h_{K_1} \leq h_{K_2}$ pointwise.

Denote by $\widetilde{\mathbb K}_n$ the collection of restrictions of the support functions 
$h_K$ to the unit sphere $S^{n-1} = \{\theta \in \R^n: |\theta| = 1\}$. A continuous function 
on $S^{n-1}$ belongs to this collection, if and only if its homogeneous extension to $\R^n$ 
represents a convex function. Thus, $\widetilde{\mathbb K}_n$ is a subset of the space 
$C(S^{n-1})$ of all continuous functions on $S^{n-1}$ which we equip with the uniform metric
$$
\|f - g\| \, = \max_{\theta \in S^{n-1}} |f(\theta) - g(\theta)|.
$$
Moreover, it follows from the definition that, for all $K_1,K_2 \in \mathbb K_n$,
$$
\rho(K_1,K_2) \, = \max_{\theta \in S^{n-1}} |h_{K_1}(\theta) - h_{K_2}(\theta)| \, = \, 
\|h_{K_1} - h_{K_2}\|.
$$
It should also be clear that if a sequence $(f_l)_{l \geq 1}$ in $\widetilde{\mathbb K}_n$ 
converges uniformly to some function $f$ as $l \rightarrow \infty$, then necessarily
$f$ belongs to $\widetilde{\mathbb K}_n$. Hence, $\widetilde{\mathbb K}_n$ is a closed 
subspace of $C(S^{n-1})$, and $({\mathbb K}_n,\rho)$ is isometric to $\widetilde{\mathbb K}_n$ 
with uniform distance.

The last inclusion in (3.2) is equivalent to the pointwise bound $h_{T(x)} \leq h_{T(x_0)} + \ep$. 
Thus, $T$ is upper semi-continuous, if and only if, for any $x_0 \in \Omega$ and $\theta \in \R^n$,
$$
\limsup_{x \rightarrow x_0}\, h_{T(x)}(\theta) = h_{T(x_0)}(\theta).
$$
This means that the real-valued function $x \rightarrow  h_{T(x)}(\theta)$ is upper semi-continuous
on $\Omega$ for any $\theta \in \R^n$, or equivalently, for any $\theta \in S^{n-1}$. Using 
the identification of ${\mathbb K}_n$ and $\widetilde{\mathbb K}_n$, we may consider a map 
$\widetilde T: \Omega \rightarrow \widetilde{\mathbb K}_n$ and say that it is upper semi-continuous, 
if for any $\theta \in S^{n-1}$,
\be
\limsup_{x \rightarrow x_0}\, \widetilde T(x)(\theta) = \widetilde T(x_0)(\theta).
\en

We need to show that any such map $\widetilde T$ is Borel measurable.
To this aim, let us recall 
that the Borel $\sigma$-algebra in $C(S^{n-1})$ is generated by the cylindrical open sets
$$
C = \big\{f \in C(S^{n-1}): f(\theta_i) < c_i \ \ i = 1,\dots,n\big\}, \quad \theta_i \in S^{n-1}, \ c_i \in \R.
$$
Assuming that (3.3) is fulfilled, the pre-image of this set
$$
 \widetilde T^{-1}(C) =  \widetilde T^{-1}(C \cap \widetilde{\mathbb K}_n) = 
\big\{x \in \Omega:  \widetilde T(x) \in C\big\}
$$
is obviously open in $\Omega$. Hence, the pre-image of any Borel set in $\widetilde{\mathbb K}_n$
is Borel in $\Omega$.
\qed

\vskip5mm
{\bf Proof of Proposition 3.1.} We apply Lemma 3.2 to the map $T(x) = \partial L(x)$
defined on the Euclidean space $\Omega = \R^n$. It was shown in \cite{R}, Corollary 24.5.1, 
that this map satisfies the property (3.2), that is, it is upper semi-continuous. Hence, it is
Borel measurable. On the other hand, it is known that the barycenter map
$$
b(K) = \int_K z\,dm_K(z), \quad K \in {\mathbb K}_n,
$$
where $m_K$ denotes a uniform distribution on $K$, is continuous with respect to the Hausdorff 
distance. Hence, the superposition $y(x) = b(\partial L(x))$ defined in (3.1) is Borel measurable 
on $\R^n$.
\qed

\vskip7mm
\section{{\bf Remarks on the Construction of Measurable Subgradients}}
\setcounter{equation}{0}

\vskip2mm
\noindent
Here we provide more remarks on the following question of independent interest which is discussed 
in the literature on non-smooth optimization.

\vskip3mm
{\bf Problem.} Given a convex function $L$ on $\R^n$, how can one construct a subgradient
$y = y(x)$ at a given point $x \in \R^n$?

\vskip3mm
As mentioned in \cite{K-Y} with reference to some recent and old publications, 
evaluating a subgradient may be a challenging task; this difficulty 
has motivated the development of numerous subdifferential approximations.
In Proposition 3.1 one  of those constructions is described using formula (3.1).
However, a preferable approach would involve directional derivatives $L(x,\theta)$ defined in (2.3) 
(which are treated as more accessible quantities).
The authors of \cite{K-Y} consider and solve this problem in dimension $n=2$.

When $L$ is differentiable everywhere, there is a unique map
$y(x) = \nabla L(x)$, and it is continuous, as was mentioned before. In dimension $n=1$, one may take
the right derivative
$$
y(x) = L'(x+) = \lim_{\ep \downarrow 0} \frac{L(x + \ep) - L(x)}{\ep}, \quad x \in \R.
$$
This number represents the right endpoint of the closed segment $\partial L(x)$. Alternatively,
one may take the left derivative of $L$ at $x$, or even more naturally -- the average of the left 
and right derivatives. It is natural to extend this approach to the high dimensional setting 
using the directional derivatives $L'(x,\theta)$ defined in (2.3). 

First note that, for $x \in E$, the gradient function is representable in terms of usual partial derivatives as
$$
\nabla L(x) = \sum_{i=1}^n \partial_i L(x)\, e_i.
$$
More generally, if $v_1,\dots,v_n$ is an orthonormal basis of $\R^n$, we have an equivalent
formula
$$
\nabla L(x) = \sum_{i=1}^n L'(x,v_i)\, v_i,
$$
which makes also sense for $x \notin E$, using the definition (2.3).
Averaging this equality over all $(v_1,\dots,v_n)$, that is, with respect to the Haar probability 
measure on the Stiefel manifolds, we obtain a Borel measurable map
\be
y(x) = n \int_{S^{n-1}} L'(x,\theta)\,\theta\,d\sigma_{n-1}(\theta), \quad x \in \R^n,
\en
where $\sigma_{n-1}$ denotes the uniform distribution on the unit sphere $S^{n-1}$. 

Since $y(x) = \nabla L(x)$ for $x \in E$, (4.1) may provide a natural generalization 
of the gradient of $L$ at a given point to the case where it does not exist in the usual sense.

\vskip3mm
{\bf Conjecture 4.1.} {\sl For the map defined in $(4.1)$ necessarily $y(x) \in \partial L(x)$ 
for all $x \in \R^n$. Is that map identical to the map defined  in $(3.1)$?
}

\vskip3mm
In dimension $n=1$, $S^0 = \{-1,1\}$ on which the uniform distribution is just the Bernoulli measure
$\sigma_0 = \frac{1}{2}\,\delta_{-1} + \frac{1}{2}\,\delta_1$. In this case, (4.1) leads to 
$$
y(x) = \frac{1}{2}\,(L'(x+) + L'(x-)).
$$

In a closely related approach consider the maps
\be
T_\ep(x) = \frac{1}{|B(x,\ep)|} \int_{B(x,\ep)} \nabla L(z)\,dz,
\en
where $B(x,\ep)$ denotes the ball in $\R^n$ with center at $x \in \R^n$ and radius $\ep>0$, with
volume $|B(x,\ep)| = \omega_n \ep^n$. The above integral is well-defined and finite. 
Here the integration may be restricted to the set $E$ of all points of differentiability of $L$, 
on which we recall that $\nabla L$ is continuous and bounded 
(due to the local Lipschitz property of $L$). 

\vskip5mm
{\bf Proposition 4.2.} {\sl  Suppose that, for every $x \in \R^n$, there exists the following limit
\be
\lim_{\ep \rightarrow 0} T_\ep(x) = y(x).
\en
Then $y$ provides a Borel measurable map such that $y(x) \in \partial L(x)$ for any $x \in \R^n$. 
}

\vskip5mm
{\bf Proof.} This is a consequence of the upper semi-continuity property of the mapping
$x \rightarrow \partial L(x)$ which we discussed before. Note that, since the map 
$x \rightarrow  \int_{B(x,\ep)} \nabla L(z)\,dz$ is continuous, $T_\ep$ is continuous as well. 
Hence, if the limit in (4.3) exists, we would obtain a Borel measurable map 
$y:\R^n \rightarrow \R^n$.

For another proof of the claim that $y(x) \in \partial L(x)$, fix $x$ and denote by $m_\ep$ 
the uniform distribution on $B(x,\ep)$ so that
\be
T_\ep(x) = \int_E \nabla L(z)\,dm_\ep(z).
\en
Note that $m_\ep(E)=1$ and recall that on the set $E$ 
$$
\left<z,\nabla L(z)\right> = L(z) + L^*(\nabla L(z)), \quad z \in E.
$$
Hence, by Jensen's inequality,
\bee
L^*(T_\ep(x)) 
 & \leq &
\int_E L^*(\nabla L(z))\,dm_\ep(z) \\
 & = & 
\int_E \big(\left<z,\nabla L(z)\right> - L(z)\big)\,dm_\ep(z) \\
 & = &
\int_E \big(\left<x,\nabla L(z)\right> - L(x)\big)\,dm_\ep(z) \\
 & & + \
\int_E \big(\left<z-x,\nabla L(z)\right> - (L(z)-L(x)\big)\,dm_\ep(z).
\ene
Here the last integral tends to zero as $\ep \rightarrow 0$, since $|z-x| < \ep$ on the ball $B(x,\ep)$
and $L$ is Lipschitz. In addition, by (4.4),
$$
\int_E \left<x,\nabla L(z)\right> dm_\ep(z) = \left<x,T_\ep(x)\right>.
$$
Thus,
$$
L^*(T_\ep(x)) \, \leq \, \left<x,T_\ep(x)\right> - L(x) + O(\ep) \quad {\rm as} \ \ep \rightarrow 0.
$$
Letting $\ep \rightarrow 0$ and using the hypothesis (4.3), we arrive at
$$
L^*(y) \, \leq \, \left<x,y\right> - L(y),  \quad y = y(x).
$$
But, by the very definition of the dual transform, the opposite inequality holds true. Hence 
the above is equivalent to
$$
L^*(y) \, = \, \left<x,y\right> - L(y),  \quad y = y(x).
$$
Finally, by Proposition 2.1, this means that $y \in \partial L(x)$.
\qed

\vskip5mm
{\bf Remark 4.3.} The formula (4.2) may be rewritten as
$$
T_\ep(x) =  \frac{n}{|B(x,\ep)|} \int_{S^{n-1}} 
\int_{B(x,\ep)} L'(z,\theta)\, \theta\,\sigma_{n-1}(\theta)\, dz.
$$
Letting $\ep \rightarrow 0$, it is natural to expect that we obtain in the limit the formula (4.1).

\vskip7mm
\section{{\bf The $\Delta_2$-Condition}}
\setcounter{equation}{0}

\vskip2mm
\noindent
Let $L$ be a non-negative convex function on $\R^n$ such that $L(0)=0$ 
and $L(x)>0$ for $x \neq 0$. It is said to satisfy the $\Delta_2$-condition, if
\be
L(2x) \leq CL(x)
\en
for all $x \in \R^n$ such that $|x| > x_0$ with some positive constants $C$ and $x_0$. Here,
choosing a larger constant $C$, one may assume without loss of generality that $x_0=1$.

On the real line a detailed treatment of this property was carried out
in the classical monograph \cite{K-R} by Krasnosel'skii and Rutickii, who considered
the so-called $N$-functions, that is, even convex functions $L:\R \rightarrow [0,\infty)$
such that $L(0) = 0$, $L(x)>0$ for $x \neq 0$, satisfying 
$$
\lim_{x \rightarrow 0} \frac{L(x)}{x} = 0 \quad {\rm and} \quad
\lim_{x \rightarrow \infty} \frac{L(x)}{x} = \infty.
$$
It is well-known that the $\Delta_2$-condition is equivalent to the assertion that
\be
\sup_{x > x_0} \, \frac{L'(x) x }{L(x)} < \infty
\en
for some $x_0>0$, where $L'(x)$ denotes the right  derivative 
of $L$ at the point $x$. In this section, we extend this characterization to the 
multidimensional setting in terms of directional derivatives.

First, let us note that the $\Delta_2$-condition (5.1) is equivalent to the property that the function
\be
R_L(r) \, = \sup_{|x| \geq 1} \, \frac{L(rx)}{L(x)},\quad r \geq 0,
\en
is finite on the positive half-axis. Let us look at the basic properties of $R_L$.

\vskip5mm
{\bf Proposition 5.1.} {\sl  Under the $\Delta_2$-condition the following properties 
hold true for the function $R = R_L$:

\vskip2mm
$a)$ \ $R$ is non-negative, non-decreasing, convex, and satisfies $R(0)=0$, $R(r)>0$ for $r > 0$;

$b)$ \ $R(1)=1$ and $1 \leq R'(1-) \leq R'(1+)$;

$c)$ \ $R(rs) \leq R(r) R(s)$ for all $r \geq 0$, $s \geq 1$;

$d)$ \ $R$ has a sub-polynomial growth: putting $p = R'(1+)$, we have
\be
R(r) \leq r^p \ (r \geq 1).
\en
}

\vskip0mm
{\bf Proof.} $a)-b)$ are clear. For $c)$, given $x \in \R^n$ with $|x| \geq 1$ and 
$r \geq 0$, $s \geq 1$, we have
$$
\frac{L(rsx)}{L(x)} = \frac{L(r(sx))}{L(sx)} \frac{L(sx)}{L(x)} \leq R(r) R(s),
$$
where we used $|sx| \geq 1$.
Taking the supremum over all admissible $x$ on the left-hand side, the claim in $c)$ follows. 

This property implies that, for any real $a \geq 1$ and any integer $m \geq 0$,
\be
R(a^m) \leq R(a)^m = \exp\{m \log R(a)\}.
\en
Given $r > 1$ and $a > 1$, we apply this bound with $m = [\frac{\log r}{\log a}] + 1$, so that
$a^m \geq r$. It gives
\be
R(r) \, \leq \,
\exp\Big\{\Big(\frac{\log r}{\log a} + 1\Big) \log R(a)\Big \} \, = \,
R(a)\,\exp\Big\{\frac{ \log R(a)}{\log a} \log r\Big \}.
\en
For $a = 1 + \ep$ with $\ep \downarrow 0$, by the Taylor expansion and using $b)$, we have
$$
R(a) = 1 + R'(1+)\, \ep + o(\ep), \quad \log R(a) = R'(1+)\, \ep + o(\ep).
$$
Hence in the limit (5.6) yields (5.4).
\qed

\vskip5mm
Note that the polynomial bound in (5.4) implies 
$$
L(x) \leq C\,|x|^p, \quad |x| \geq 1,
$$
with constant $C = \sup_{|x| = 1} L(x)$. 

One can now turn to the characterization of the $\Delta_2$-condition. Recall that $L'(x,x)$ 
denotes the directional derivative of $L$ at the point $x$ in the direction of $x$, according to (2.3). 

\vskip5mm
{\bf Proposition 5.2.} {\sl The function $L$ satisfies the $\Delta_2$-condition, if and only if
\be
\sup_{|x| \geq 1} \, \frac{L'(x,x)}{L(x)} \, = \,
\sup_{|x| \geq 1} \, \sup_{y \in \partial L(x)} \frac{\left<x,y\right>}{L(x)} \, < \, \infty.
\en
In this case, both suprema are at most $p = R_L'(1+)$. If $L$ is differentiable,
this condition simplifies to
\be
\sup_{|x| \geq 1} \, \frac{\left<\nabla L(x),x\right>}{L(x)} < \infty.
\en
}

\vskip2mm
Since $\left<x,y\right> = L(x) + L^*(y)$ for $y \in \partial L(x)$, (5.7) may be equivalently 
stated as
$$
\sup_{|x| \geq 1} \, \sup_{y \in \partial L(x)} \frac{L^*(y)}{L(x)} < \infty.
$$

Clearly, in the class of $N$-functions in dimension $n=1$, (5.7)-(5.8) are reduced to the
characterization (5.2).

\vskip3mm
{\bf Proof of Proposition 5.2.} In Section 2 we mentioned a general identity
$$
L'(x,\theta) = h_{\partial L(x)}(\theta) = \sup_{y \in \partial L(x)} \left<y,\theta\right>, \quad
x,\theta \in \R^n,
$$
which for $\theta = x$ becomes
$$
L'(x,x) =  \sup_{y \in \partial L(x)} \left<y,x\right>.
$$
This explains the equality in (5.7).

Let us now explain why the finiteness ofboth suprema in (5.7) is necessary for the 
$\Delta_2$-condition. Assume that (5.1) holds. By (5.4), 
\be
L(rx) \leq r^p L(x), \quad x \in \R^n, \ |x| \geq 1, \ r \geq 1.
\en
By the convexity of $L$, for all $h \in \R^n$ and $y \in \partial L(x)$,
$$
L(x+h) - L(x) \geq \left<h,y\right>.
$$
Applying this bound with $h = \ep x$ and choosing in (5.9) $r = 1+\ep$, $\ep>0$, we get
$$
\ep \left<x,y\right> \leq \big((1+\ep)^p - 1\big) L(x).
$$
Therefore
$$
\sup_{|x| \geq 1} \, \sup_{y \in \partial L(x)} \frac{\left<x,y\right>}{L(x)} \leq
\frac{1}{\ep}\,\big((1+\ep)^p - 1\big).
$$
Letting $\ep \rightarrow 0$, we arrive at 
$$
\sup_{x \neq 0} \, \sup_{y \in \partial L(x)} \frac{\left<x,y\right>}{L(x)} \leq p,
$$
and the finiteness of the second supremum in (5.7) follows. 
Similarly, from (5.9) with $r = 1+\ep$ it follows that
$$
L'(x,x) \, = \, \lim_{\ep \downarrow 0} \frac{L(x + \ep x) - L(x)}{\ep} \, \leq \,
\lim_{\ep \downarrow 0} \frac{(1+\ep)^p L(x) - L(x)}{\ep} \, = \, p L(x).
$$
That is, we obtain (5.7) and also notice that the supremum in it is at most $p$.

Now, let us start with the hypothesis (5.7) and denote by $c$ the first supremum.
Fix $x\in \R^n$ such that $|x| \geq 1$ and 
consider the function $V(r) = L(rx)$, $r \geq 1$. It is convex, non-decreasing, and has 
the right derivative
\bee
V'(r+) 
 & = &
\lim_{\ep \downarrow 0} \frac{V(r + \ep) - V(r)}{\ep} \\
 & = &
\lim_{\ep \downarrow 0} \frac{L(rx + \ep x) - L(rx)}{\ep} \\
 & = & 
L'(rx,x) \ = \ \frac{L'(rx,rx)}{r} \ \leq \ \frac{c}{r}\, L(rx).
\ene
Hence, $V$ satisfies a differential inequality $V'(r+) \leq \frac{c}{r}\,  V(r)$. 
Since $V(r) > 0$, this inequality is equivalent to the property that the function 
$\log V(r) - c \log r$ is non-increasing in $r \geq 1$. In particular,
$V(2) \leq 2^c\, V(1)$, which amounts to the definition (5.1) with $C = 2^c$ and $x_0=1$.
\qed

\vskip7mm
\section{{\bf  A Two-Sided $\Delta_2$-Condition}}
\setcounter{equation}{0}

\vskip2mm
\noindent
Again, let $L$ be a non-negative convex function on $\R^n$ such that $L(0)=0$, 
$L(x)>0$ for $x \neq 0$. In order to describe the possible behavior of $L(x)$ for large $x$, 
Krasnosel'skii and Rutickii also discussed in \cite{K-R} related properties such as 
the $\Delta'$-condition and the $\Delta_3$-con\-dition.
Here we focus on a formally different class of convex functions, whose behaviour can be
controlled not only at infinity, but also around zero.

\vskip5mm
{\bf Definition 6.1.} We say that $L$ satisfies a two-sided $\Delta_2$-condition, if
\be
L(2x) \leq CL(x) \quad {\rm for \ all} \ x \in \R^n
\en
with some constant $C$.

\vskip3mm 
It was shown in \cite{K-R} that, if $L$ satisfies (6.1) for all $x$ large enough (in dimension $n=1$),
it is possible to redefine $L$ near zero in such a way that this inequality would hold
for all $x$. The new function is thus equivalent to the original one (in the sense of \cite{K-R}), 
and therefore both functions lead to the identical associated Orlicz spaces. 
Nevertheless, this two-sided $\Delta_2$-condition introduces some new remarkable 
features to the analysis of convex functions.

Note that (6.1) is fulfilled if and only if the function
\be
\Phi_L(r) \, = \sup_{x \neq 0} \, \frac{L(rx)}{L(x)},\quad r \geq 0,
\en
is finite on the positive half-axis. We call it the Young function associated to $L$.
Basic properties of $\Phi_L$ are almost identical to those of $R_L$. 

\vskip5mm
{\bf Proposition 6.2.} {\sl  Under the two-sided $\Delta_2$-condition the following properties 
hold true for the function $\Phi = \Phi_L$:

\vskip2mm
$a)$ \ $\Phi$ is non-negative, non-decreasing, convex, and satisfies $R(0)=0$, $R(r)>0$ for $r > 0$;

$b)$ \ $\Phi(1)=1$ and $1 \leq \Phi'(1-) \leq \Phi'(1+)$;

$c)$ \ $\Phi$ is sub-multiplicative: $\Phi(rs) \leq \Phi(r) \Phi(s)$ for all $r,s \geq 0$;

$d)$ \ $\Phi$ has a sub-polynomial growth: putting $p_- = \Phi'(1-)$ and $p_+ =\Phi '(1+)$, we have
\be
\Phi(r) \leq r^{p_-} \ (0 \leq r\leq 1),  \qquad \Phi(r) \leq r^{p_+} \ (r \geq 1).
\en
}

\vskip3mm
Properties $a)-b)$ are obvious. The proof of $c)$ and of the second inequality in (6.3) 
is similar to the one from Proposition 5.1.
As for the first inequality, note that, due to the improved property $c)$,
the bound (5.5) holds true for all $a \geq 0$. If $0<r<1$ and $0<a<1$, we apply this bound with 
$m = [\frac{\log r}{\log a}]$, so that $a^m \geq r$. Using the monotonicity of $R$, it gives
$$
R(r) \, \leq \, R(a)^m \, \leq \,
\exp\Big\{\Big(\frac{\log r}{\log a} - 1\Big) \log R(a)\Big \} \, = \,
\frac{1}{R(a)}\,\exp\Big\{\frac{ \log R(a)}{\log a} \log r\Big \}.
$$
For $a = 1 - \ep$ with $\ep \downarrow 0$, we have
$$
R(a) = 1 - R'(1-)\, \ep + o(\ep), \quad \log R(a) = -R'(1-)\, \ep + o(\ep).
$$
Therefore, in the limit the above inequality yields the first bound in (6.3).

Now, as a closely related function, also introduce
\be
\Psi_L(r) \, = \inf_{x \neq 0} \, \frac{L(rx)}{L(x)},\quad r \geq 0.
\en
From (6.2) it follows immediately that
\be
\Psi_L(r) = \frac{1}{\Phi_L(1/r)}, \quad r>0.
\en
Let us record without proof its basic properties.

\vskip5mm
{\bf Proposition 6.3.} {\sl  Under the two-sided $\Delta_2$-condition the following 
properties hold true for the function $\Psi = \Psi_L$:

\vskip2mm
$a)$ \ $\Psi$ is non-negative, non-decreasing, and satisfies $\Psi(0)=0$, $\Psi(r)>0$ for $r > 0$;

$b)$ \ $\Psi(1)=1$ and $\Psi'(1-) = \Phi'(1+)$, $\Psi'(1+) = \Phi'(1-)$;

$c)$ \ $\Psi$ is super-multiplicative: $\Psi(rs) \geq \Psi(r) \Psi(s)$ for all $r,s \geq 0$;

$d)$ \ $\Psi$ has a super-polynomial growth: 
\be
\Psi(r) \geq r^{p_+} \ (0 \leq r\leq 1),  \qquad \Psi(r) \geq r^{p_-} \ (r \geq 1).
\en
}

\vskip2mm
The latter bounds follow from (6.5) and (6.3).

Since $\Psi(r) \leq \Phi(r)$, the comparison of (6.6) with (6.3) leads to the following conclusion.

\vskip5mm
{\bf Corollary 6.4.} {\sl  Under the two-sided $\Delta_2$-condition, the function $\Phi_L$
is differentiable at the point $r=1$ and has derivative $p = \Phi_L'(1)$, if and only if 
$$
\Phi_L(r) = \Psi_L(r) = r^p
$$ 
for all $r \geq 0$, that is, if $L$ is positive homogeneous of order $p$:
\be
L(rx) = r^p L(x), \quad x \in \R^n, \ r \geq 0.
\en
}

\vskip0mm
Necessarily $p = p_- = p_+$.
Typical examples of such functions are provided by $L(x) = \|x\|^p$, $x \in \R^n$,
with parameter $p \geq 1$, where $\|\cdot\|$ is an arbitrary norm on $\R^n$. 

Let us also mention several non-homogeneous examples of the form
$$
L(x) = V(\|x\|), \quad x \in \R^n,
$$
where $V$ is a Young function, that is, $V:[0,\infty) \rightarrow [0,\infty)$ is convex, $V(0) = 0$, 
$V(r) > 0$ for $r>0$. Then
$$
\Phi_L(r) = \Phi_V(r) = \sup_{s > 0} \, \frac{V(rs)}{V(s)},\quad r \geq 0.
$$

\vskip2mm
{\bf Example 6.5.} The functions 
$$
V(r) = r^p \max(r,1)
$$
with parameter $p \geq 1$ are convex and sub-multiplicative. Hence,
the associated Young functions do not change: $\Phi_L = \Phi_V = V$.
In this case, 
$$
p_- = V'(1-) = p \quad {\rm and} \quad p_+= V'(1+) = p+1. 
$$
Also, by the formula (6.5),
$$
\Psi_L(r) = \Psi_V(r) = r^p \min(r,1).
$$
Note that this function is not convex.

\vskip2mm
{\bf Example 6.6.} 
For the functions
$$
V(r) = r^p \log(1+r), \quad p \geq 1, 
$$
elementary comptutations show that 
\be
\Phi_L(r) = \Phi_V(r) = r^p \max(r,1)
\en
like in the previous example. For the proof, fix $r>0$ and consider the function 
$$
U(s) = \log(1 + rs) - C\log(1 + s), \quad s \geq 0,
$$ 
where 
$C>0$ is a parameter. We need to find the smallest value of $C = C_r$ such that 
$U(s) \leq 0$ for all $s \geq 0$. Since $U(0) = 0$, it should be required that $U'(0) \leq 0$. 
We have 
$$
U'(s) = \frac{r}{1+rs} - \frac{C}{1+s}
$$ 
and $U'(0) = r-C$. Hence, necessarily 
$C \geq r$. With the choice $C=r$, we have $U'(s) \leq 0$ whenever $r \geq 1$, and then
$$
\Phi_V(r) \, = r^p \sup_{s > 0} \, \frac{\log(1+rs)}{\log(1+s)} = r^{p+1},\quad r \geq 1.
$$
In the other case $r < 1$, the last ratio is smaller than 1 and tends to 1 as $s \rightarrow \infty$. 
As a result, we arrive at (6.8). 

\vskip2mm
{\bf Example 6.7.} 
In the similar example $V(r) = r^p \log(2+r)$, we have
\be
c_1 V(r) \leq \Phi_V(r) \leq c_2 V(r)
\en
with some absolute constants $c_1,c_2 > 0$. Indeed, by the definition,
$$
\Phi_V(r) = x^p\, \sup_{s > 0}\, U_r(s), \quad U_r(s) = \frac{\log(2+rs)}{\log(2+s)}.
$$
The last supremum is at least $U_r(1) = \log(2+r)/\log 3$, which yields the lower bound in (6.9)
with $c_1 = 1/\log 3$. Using $2+rs \leq (2+r)(2+s)$, we get 
$$
U_r(s) \leq 1 + \frac{\log(2+r)}{\log(2+s)} \leq 1 + \frac{\log(2+r)}{\log 2}
$$
which yields the upper bound with $c_2 = 1 + 1/\log 2$.

\vskip7mm
\section{{\bf Characterization of the Two-Sided $\Delta_2$-Conditions}}
\setcounter{equation}{0}

\vskip2mm
\noindent
Applying Proposition 6.2 in analogy with Proposition 5.2, one may similarly involve directional
derivatives and subgradients in order to give necessary and sufficient conditions
for $L$ to satisfy the two-sided $\Delta_2$-condition (6.1). 

Recall that 
$L'(x,x)$ denotes the directional derivative of $L$ in the direction of $x$ defined in $(2.3)$
and $\Phi_L$ denotes the associated Young function defined in (6.2).

\vskip5mm
{\bf Proposition 7.1.} {\sl  Let $L$ be a non-negative convex function on $\R^n$ 
such that $L(0)=0$ and $L(x)>0$ for $x \neq 0$. The function 
$L$ satisfies the two-sided $\Delta_2$-condition, if and only if
\be
\sup_{x \neq 0} \, \frac{L'(x,x)}{L(x)} \, = \,
\sup_{x \neq 0} \, \sup_{y \in \partial L(x)} \frac{\left<x,y\right>}{L(x)} \, < \, \infty.
\en
In this case, both suprema are at most $p_+ = \Phi_L'(1+)$. If $L$ is differentiable,
this condition simplifies to
\be
\sup_{x \neq 0} \, \frac{\left<\nabla L(x),x\right>}{L(x)} < \infty.
\en
}

\vskip2mm
Since $\left<x,y\right> = L(x) + L^*(y)$ for $y \in \partial L(x)$, (7.1) may be equivalently 
stated as
$$
\sup_{x \neq 0} \, \sup_{y \in \partial L(x)} \frac{L^*(y)}{L(x)} < \infty.
$$

The proof of Proposition 7.1 is similar to the proof of Proposition 5.2, so we omit it.

Another interesting question we need to address is when the Legendre transform
$$
L^*(x) = \sup_y\, \big[ \left<x,y\right> - L(y)\big]
$$
satisfies the two-sided $\Delta_2$-condition. Recall that 
$1 \leq p_- \leq p_+$ where $p_- = \Phi_L'(1-)$ and $p_+ = \Phi_L'(1+)$.
As before, we assume that $L$ is a non-negative convex function on $\R^n$ 
such that $L(0)=0$ and $L(x)>0$ for $x \neq 0$. 

\vskip5mm
{\bf Proposition 7.2.} {\sl If $L$ satisfies the two-sided $\Delta_2$-condition with $p_- > 1$, 
then $L^*$ satisfies the two-sided $\Delta_2$-condition as well. In this case
\be
\Phi_{L^*}(1+) \leq \frac{p_-}{p_- - 1}, \quad \Phi_{L^*}(1-) \leq \frac{p_+}{p_+ - 1}.
\en
}

\vskip2mm
{\bf Proof.} Suppose that $p = p_- > 1$ and let $r \geq 1$, $x \in \R^n$.
Applying the first inequality in (6.3), i.e. $L(sx) \leq s^p L(x)$, $s \in [0,1]$, we get
\bee
L^*(rx) 
 & = &
\sup_y\, \big[ \left<rx,y\right> - L(y)\big] \, = \,
\sup_y\, \big[ \left<x,y\right> - L(y/r)\big] \\
 & \geq &
\sup_y\, \big[ \left<x,y\right> - r^{-p}L(y)\big] \\
 & = &
r^{-p}\, \sup_y\, \big[ \left<r^p x,y\right> - L(y)\big] \, = \, r^{-p}\,L^*(r^p x).
\ene
Thus,
$$
L^*(r^p x) \leq r^p L^*(rx).
$$
Changing the variable $rx = y$, we get $L^*(r^{p-1} y) \leq r^p L^*(y)$. The replacement $r^{p-1} = s$
yields $L^*(s y) \leq s^{q_+} L^*(y)$ or equivalently
$$
L^*(r x) \leq r^{q_+} L^*(x), \quad r \geq 1, \ q_+ = \frac{p_-}{p_--1}.
$$
This implies that $L^*$ satisfies the two-sided $\Delta_2$-condition and also proves the first inequality 
in (7.3).

A similar argument based on the second inequality in (6.3) leads to 
$$
L^*(r x) \leq r^{q_-} L^*(x), \quad 0 \leq r \leq 1, \ q_- = \frac{p_+}{p_+-1},
$$
which gives the second inequality in (7.3).
\qed

\vskip5mm
{\bf Remark 7.3.} If $L(x) \geq c|x|$ for all $x \in \R^n$ with some constant $c>0$, 
then $L^*(x) = 0$ for all sufficiently small $|x|$. Hence, in this case we may not speak 
about the $\Delta_2$-condition for the Legendre transform of $L$. For an illustration, 
let us return to Example 6.5 with $p=1$ and consider the convex function
$$
L(x) = |x|\,\max(|x|,1), \quad x \in \R^n.
$$
It satisfies the $\Delta_2$-condition with $\Phi_L(r) = r\,\max(r,1)$, in which case $p_- = 1$, $p_+ = 2$.
Then, as easy to see,
$$
L^*(x) = \max\big\{0, |x|-1, |x|^2/4\big\}.
$$
This function is vanishing for $|x| \leq 1$.

Finally, let us refine Proposition 7.2 for the class of homogeneous convex functions.

\vskip5mm
{\bf Proposition 7.4.} {\sl If $L$ is positive homogeneous of order $p \geq 1$, it
satisfies the two-sided $\Delta_2$-condition. In this case, $L^*$ satisfies the two-sided 
$\Delta_2$-condition
if and only if $p > 1$, and then it is positive homogeneous of order $q = p/(p-1)$.
}

\vskip5mm
{\bf Proof.} Let $L(rx) = r^p L(x)$ for all $r \geq 0$ and $x \in \R^n$.
As in the previous proof, if $r>0$, 
\bee
L^*(rx) 
 & = &
\sup_y\, \big[ \left<rx,y\right> - L(y)\big] \, = \,
\sup_y\, \big[ \left<x,y\right> - L(y/r)\big] \\
 & = &
\sup_y\, \big[ \left<x,y\right> - r^{-p}L(y)\big] \\
 & = &
r^{-p}\, \sup_y\, \big[ \left<r^p x,y\right> - L(y)\big] \, = \, r^{-p}\,L^*(r^p x).
\ene
Thus,
$$
L^*(r^p x) = r^p L^*(rx).
$$
If $p>1$, this means that $L^*$ is positive homogeneous of order $q = p/(p-1)$.
If $p=1$, the above equality is only possible when $L^*(x) = 0$ or $L^*(x) = \infty$
for all $x \in \R^n$.
\qed

\vskip7mm
\section{{\bf Luxemburg and Orlicz Pseudo-Norms for Vector-Valued Functions}}
\setcounter{equation}{0}

\vskip2mm
\noindent
Let $L:\R^n \rightarrow [0,\infty]$ be a lower semi-continuous convex function such that 
$L(0)=0$, $L(x)>0$ for all $x \neq 0$, which is finite in some neighborhood of the origin.
The conjugate transform
$$
L^*(x) = \sup_{y \in \R^n} \big[\left<x,y\right> - L(y)\big], \quad x \in \R^n,
$$
is also non-negative, convex, but may vanish near zero and take an infinite value.
The following definition extends the classiocal notions of Luxemburg and Orlicz norms 
from the case of scalar functions to the multidimensional setting.

\vskip5mm
{\bf Definition 8.1.} Let $(\Omega,\frak M,\lambda)$ be a probability space. Given a measurable 
function $u:\Omega \rightarrow \R^n$, define its Luxemburg pseudo-norm by
\be
\|u\|_L = \|u\|_{L(\lambda)} = \inf\Big\{r > 0: \int L(u/r)\,d\lambda \leq 1\Big\},
\en
when the value of $r$ under the infimum sign exists, and $\|u\|_L = \infty$ otherwise.

If $\|u\|_L$ is finite, define the Orlicz pseudo-norm of $u$ by
\be
|u|_L = |u|_{L(\lambda)} = \sup\Big\{\int \left<u,v\right>d\lambda: \|v\|_{L^*} \leq 1\Big\}.
\en

\vskip3mm
The scalar case in  (8.1)--(8.2) was discussed in detail in Krasnoselskii and Rutickii \cite{K-R} 
and Maligranda \cite{M}. The latter book also contains historical remarks: The definition (8.2) 
was first given in 1930's by Orlicz \cite{O1} assuming the $\Delta_2$-condition and then 
in \cite{O2} for general Young functions. The equivalent definition (8.1) was given later in 1950's 
in the work by Nakano \cite{N} and Luxemburg \cite{Lu}. 
In the literature, the use of (8.2) is rather rare, and the equality
(8.1) is often taken as the definition of the Orlicz norm.

In (8.2), $\|v\|_{L^*}$ is defined according to (8.1), although the condition 
$L^*(x)>0$ for all $x \neq 0$ may be violated (cf. Remark 7.3).
If $\|u\|_L$ is finite and positive, the infimum in (8.1) is attained at $r = \|u\|_L$, by applying
the dominated  convergence theorem and using the property that $r \rightarrow L(rx)$ 
is non-decreasing in $r \geq 0$. Thus
$$
\int L(u/\|u\|_L)\,d\lambda = 1 \quad (0 < \|u\|_L < \infty).
$$

Note that $\|u\|_L = 0$ implies that $u=0$ with $\lambda$-probability 1. Indeed, suppose
that $u \neq 0$ on a set $A \subset \Omega$ of positive $\lambda$-probability. Then, using
$L(x)>0$ for all $x \neq 0$, we would have $L(\alpha u) \uparrow \infty$ as $\alpha \uparrow \infty$
on the set $A$, so that $\int L(\alpha u)\,d\lambda \uparrow \infty$ or equivalently
$\int L(u/r)\,d\lambda \uparrow \infty$ as $r \downarrow 0$. According to (8.1), we then get
$\|u\|_L > 0$.

As will be clarified later on, the integrals in (8.2) are well-defined and bounded
from above when $\|u\|_L$ is finite. This restriction may be removed, since in this definition, 
one may additionally require that under the supremum sign the functions $v$ satisfy 
$\left<u,v\right> \geq 0$. Indeed, put $\widetilde v(\omega) = v(\omega)$ if
$\left<u(\omega),v(\omega)\right> \geq 0$ and $\widetilde v(\omega) = 0$ otherwise. Then
$$
\int \left<u,v\right>d\lambda \leq \int \left<u,\widetilde v\right> d\lambda, \quad 
L^*(\widetilde v/r) \leq L^*(v/r)
$$
for any $r>0$, and hence $\|\widetilde v\|_{L^*} \leq \|v\|_{L^*}$. So, after this modification
of $v$ the integrals in (8.2) may only increase, while the norm of the modified function
$\widetilde v$ may only decrease. Therefore,
\begin{eqnarray}
|u|_L
 & = &
\sup\Big\{\int \left<u,v\right>d\lambda: \|v\|_{L^*} \leq 1,\, \left<u,v\right> \geq 0\Big\} \nonumber \\
 & = &
\sup\Big\{\int \left<u,v\right>^+ d\lambda: \|v\|_{L^*} \leq 1\Big\},
\end{eqnarray}
which may be taken as an equivalent definition instead of (8.2). With this approach there is
no need to pose the restriction that $\|u\|_L$ is finite (although, as we will see,
the finiteness of $|u|_L$ is equivalent to the finiteness of $\|u\|_L$).

Note that $|u|_L = 0$ if and only if  $\left<u,v\right> \leq 0$ $\lambda$-a.e. whenever  
$\|v\|_{L^*} \leq 1$. By homogeneity of this pseudo-norm, we get $\left<u,v\right> \leq 0$ $\lambda$-a.e.
for any bounded $v$, which is only possible when $u=0$ $\lambda$-a.e., similarly
to the Luxemburg pseudo-norm.

Let us now list basic properties of the functionals in (8.1)-(8.2).

\vskip5mm
{\bf Proposition 8.2.} {\sl 
The functional $\|u\|_L$ is non-negative, positive homogeneous of order $1$ and convex:

\vskip2mm
$a)$ \ $0 \leq \|u\|_L \leq \infty$;

\vskip1mm
$b)$ \ $\|\alpha u\|_L = \alpha\, \|u\|_L$ for all $\alpha>0$;

\vskip1mm
$c)$ \ $\|t_1 u_1 + t_2 u_2\|_L  \leq t_1 \|u_1\|_L + t_2 \|u_2\|_L$ for all $t_1,t_2 \geq 0$ 
such that $t_1+t_2=1$.

\vskip1mm
$d)$\, Hence, this functional is subadditive: $\|u_1 + u_2\|_L \leq \|u_1\|_L + \|u_2\|_L$.

\vskip2mm
\noindent
The same properties are fulfilled for the functional $|u|_L$.
}

\vskip5mm
{\bf Proof.} $a)$ is clear. To show $b)$, note that, by the definition,
\bee
\|\alpha u\|_L
 & = &
\inf\Big\{r > 0: \int L(\alpha u/r)\,d\lambda \leq 1\Big\} \\
 & = &
\inf\Big\{\alpha r' > 0: \int L(u/r')\,d\lambda \leq 1\Big\} \, = \, \alpha\, \|u\|_L.
\ene
To prove $c)$, we may assume that $t_1,t_2>0$ and that the pseudo-norms
$\|u_1\|_L = r_1$ and $\|u_2\|_L = r_2$
are positive and finite. By the convexity of $L$,
\bee
L\Big(\frac{t_1 u_1 + t_2 u_2}{t_1 r_1 + t_2 r_2}\Big)
 & = &
L\Big(\frac{t_1 r_1\,\frac{u_1}{r_2} + t_2 r_2\,\frac{u_2}{r_2}}{t_1 r_1 + t_2 r_2}\Big) \\
 & \leq &
\frac{t_1 r_1}{t_1 r_1 + t_2 r_2}\,L(u_1/r_1) + \frac{t_2 r_2}{t_1 r_1 + t_2 r_2}\,L(u_2/r_2).
\ene
Integrating both sides over $\lambda$ and using $\int L(u_1/r_1)\,d\lambda = \int L(u_2/r_2)\,d\lambda = 1$,
we get
$$
\int L\Big(\frac{t_1 u_1 + t_2 u_2}{t_1 r_1 + t_2 r_2}\Big)\,d\lambda \leq 1.
$$
This is the required inequality in $c)$. 
In view of the homogeneity, we have the subadditivity as well.
The claim about $|u|_L$ is similar.
\qed

\vskip5mm
Note that since $L$ is not required to be even, i.e. $L(-x) = L(x)$,
$x \in \R^n$, the property $\|-u\|_L = \|u\|_L$ may not hold. 
However, if $L$ is even, then we obtain a norm in the Orlicz Banach space of all measurable 
functions $u$ on $\Omega$ such that $\|u\|_L < \infty$.

\vskip5mm
{\bf Example 8.3.} Restricting ourselves to dimension $n=1$, the most frequent choice of the convex 
function $L$ is the power function $L(x) = |x|^p$, $1 < p < \infty$. Then the dual transform is given
by a multiple of a power function, namely
$$
L^*(x) = \frac{1}{q p^{q-1}}\, |x|^q, \quad q = \frac{p}{p-1}.
$$
Hence
$$
\|u\|_L = \|u\|_p = \Big(\int |u|^p\,d\lambda\Big)^{1/p}, \quad 
\|u\|_{L^*} = \frac{1}{p^{1/p}\, q^{1/q}}\,\|u\|_q,
$$
and
\be
|u|_L = p^{1/p}\, q^{1/q}\,\|u\|_p = p^{1/p}\, q^{1/q}\,\|u\|_L.
\en

\vskip7mm
\section{{\bf Relationship Between Luxemburg and Orlicz Pseudo-Norms}}
\setcounter{equation}{0}

\vskip2mm
\noindent
First let us emphasize the following natural bound on the integrals in (8.2).
We keep the same assumptions about the convex function $L$ as in the previous section.
As is standard, put $a^+ = \max(a,0)$ for $a \in \R$.

\vskip5mm
{\bf Proposition 9.1.} {\sl Given two measurable functions $u,v:\Omega \rightarrow \R^n$,
we have
\be
\int \left<u,v\right>^+  d\lambda \leq 2\,\|u\|_L \, \|v\|_{L^*}.
\en
}

\vskip2mm
{\bf Proof.} The cases where one or both of the pseudo-norms $\|u\|_L$ and $\|v\|_{L^*}$ are zero
or infinite are obvious. Hence, by homogeneity of (9.1) with respect to $u$ and $v$, we may assume 
that $\|u\|_L = \|v\|_{L^*} = 1$. Recall that
\be
\left<x,y\right> \leq L(x) + L^*(y), \quad x,y \in \R^n.
\en
Since the right-hand side is non-negative, this inequality may be sharpened to 
$\left<x,y\right>^+ \leq L(x) + L^*(y)$. In particular, 
$$\left<u,v\right>^+ \leq L(u) + L^*(v).
$$
Integrating this inequality over $\lambda$, we arrive at (9.1).
\qed

\vskip5mm
{\bf Proposition 9.2.} {\sl Suppose additionally that $L$ is everywhere finite.
For any measurable function $u:\Omega \rightarrow \R^n$, 
\be
\|u\|_L \leq |u|_L \leq 2\,\|u\|_L.
\en
}

\vskip4mm
{\bf Proof.}
By Proposition 9.1, if $\|v\|_{L^*} \leq 1$, then
$$
\int \left<u,v\right>^+  d\lambda \leq 2\,\|u\|_L.
$$
According to (8.3), this yields $|u|_L \leq 2\,\|u\|_L$ which is the second inequality in (9.3). 

To derive the first inequality, it is sufficient to consider the case where $|u|_L$ is finite and positive.
Hence, we need to show that
\be
\int L(u/|u|_L)\,d\lambda \leq 1.
\en

{\bf Step 1}: Reduction of (9.4) to bounded functions $u$.
Consider the sequence of functions $u_k$ defined to be $u$ 
if $|u| \leq k$ and to be zero otherwise. Putting $\Omega_k = \{|u| \leq k\}$, we have, by (8.3),
$$
|u_k|_L = \sup\Big\{\int_{\Omega_k} 
\left<u,v\right>d\lambda: \|v\|_{L^*} \leq 1,\, \left<u,v\right> \geq 0\Big\},
$$
which readily implies $|u_k|_L \leq |u|_L$. As a consequence,
$L(u/|u|_L) \leq L(u_k/|u_k|_L)$ on $\Omega_k$.
Integrating this inequality over the measure $\lambda$ and applying (9.4) to $u_k$, we get
$$
\int_{\Omega_k} L\Big(\frac{u}{|u|_L}\Big)\,d\lambda  \leq 
\int_{\Omega_k} L\Big(\frac{u_k}{|u_k|_L}\Big)\,d\lambda  = 
\int_\Omega L\Big(\frac{u_k}{|u_k|_L}\Big)\,d\lambda \leq 1.
$$
It remains to send $k \rightarrow \infty$, and then we obtain the desired relation (9.4) for
the function $u$.

\vskip3mm
{\bf Step 2}: The case of bounded $u$. By homogeneity, we may assume that $|u|_L = 1$, 
and then (9.4) is reduced to
\be
\int L(u)\,d\lambda \leq 1.
\en

Recall that, by Proposition 2.1,
\be
\left<x,y\right> = L(x) + L^*(y), \quad x \in \R^n, \ y \in\partial L(x).
\en
Consider the function $v = y(u)$, where $y$ is a Borel measurable map from Proposition 3.1
(necessarily $v = \nabla L(u)$ when $L$ is differentiable at $u$). This is the only place where the Borel
measurability of $y$ is needed in order to ensure that $v$ is measurable on $\Omega$. 
In addition, applying the upper bounds (2.2) and (2.6) with $r=1$,
we get
$$
v \, \leq \, \sup_{|h| \leq 1} \big(L(u+h) - L(u)\big)
$$
and
$$
L^*(v) \, \leq \, |u|\, \sup_{|h| \leq 1} \big(L(u+h) - L(u)\big),
$$
which implies that both $v$ and $L^*(v)$ are bounded as well.

Applying the equality (9.6) to the couple $(u,v)$ in place
of $(x,y)$ and integrating over $\lambda$, it follows that
\be
\int \left<u,v\right> d\lambda = \int L(u)\,d\lambda + \int L^*(v)\,d\lambda.
\en
Note that all integrands are bounded non-negative functions, so that the integrals are finite
and non-negative.

If $\int L^*(v)\,d\lambda \leq 1$, that is, when $\|v\|_{L^*} \leq 1$, then, by the definition (8.2), 
$$
\int \left<u,v\right> d\lambda \leq |u|_L = 1.
$$  
In view of (9.7),  this implies that $\int L(u)\,d\lambda \leq 1$, proving (9.5) in this case.

In the other case where $c = \int L^*(v)\,d\lambda > 1$, consider the function $\widetilde v = v/c$.
Since $L^*$ is convex and $L^*(0)=0$, we have $L^*(\alpha x) \leq \alpha L^*(x)$ for all $x \in \R^n$
and $0 \leq \alpha \leq 1$. In particular,
$$
\int L^*(\widetilde v)\,d\lambda \leq \frac{1}{c} \int L^*(v)\,d\lambda = 1.
$$
Then again by (8.2), $\int \left<u,\widetilde v\right> d\lambda \leq |u|_L = 1$, that is,
$$
\int \left<u,v\right> d\lambda \leq \int L^*(v)\,d\lambda.
$$ 
By (9.7), this yields $\int L(u)\,d\lambda \leq 0$, which is not possible. Thus (9.5) is proved.
\qed

\vskip5mm
{\bf Remark 9.3.} By the one-dimensional Example 8.3, the factor 2 cannot be improved
in (9.3). Indeed, the constant $p^{1/p}\, q^{1/q}$ in the equality (8.4) is greater than 1
and takes the maximal value 2 for $p=q=2$.

\vskip7mm
\section{{\bf Perturbations of the Luxemburg Pseudo-Norm}}
\setcounter{equation}{0}

\vskip2mm
\noindent 
In various problems based, for example, on the application of Hamilton-Jacobi equations,
one has to require a super-linear growth condition for $L$:
$$
\frac{L(x)}{|x|} \rightarrow \infty \quad {\rm as} \ |x| \rightarrow \infty.
$$
However, this condition might be unnecessary for transport-energy bounds. In order to reduce
such inequalities to the case where the convex cost functions satisfy this growth condition,
one may approximate $L$ with convex functions 
$$
L_\ep(x) = L(x) + \ep L_0(x), \quad x \in \R^n, \ \ep>0,
$$
where $L_0$ satisfies the required growth condition, for example, $L_0(x) = \frac{1}{2}\,|x|^2$.

Thus, let $L:\R^n \rightarrow [0,\infty]$ be a convex function which is finite in a neigborhood 
of zero and such that $L(0)=0$ and $L(x)>0$ for all $x \neq 0$, and let
$L_0:\R^n \rightarrow [0,\infty)$ be a convex function such that $L_0(0)=0$ 
and $L_0(x)>0$ for all $x \neq 0$. Here we collect some useful properties of perturbed convex functions.

\vskip5mm
{\bf Proposition 10.1.} {\sl Suppose that $L_0$ satisfies the super-linear growth condition.
The following properties hold true as $\ep \downarrow 0$:

\vskip3mm
$a)$ \ $L(x) \leq L_\ep(x)$ and $L_\ep(x) \downarrow L(x)$ for all $x \in \R^n$;

\vskip2mm
$b)$ \ $L_\ep^*(x) \leq L^*(x)$ and $L_\ep^*(x) \uparrow L^*(x)$ for all $x \in \R^n$;

\vskip2mm
$c)$ \ Given a measurable function $u:\Omega \rightarrow \R^n$, we have
$\|u\|_{L(\lambda)} \leq \|u\|_{L_\ep(\lambda)}$ for all $\ep>0$, and
if $u$ is bounded, then $\|u\|_{L_\ep(\lambda)} \downarrow \|u\|_{L(\lambda)}$.

\vskip2mm
$d)$ \ Given a measurable function $u:\Omega \rightarrow \R^n$, we have
$$ 
\|u\|_{L_\ep^*(\lambda)} \leq \|u\|_{L^*(\lambda)} \quad {\sl and} \quad
\|u\|_{L_\ep^*(\lambda)} \uparrow \|u\|_{L^*(\lambda)}.
$$
}

\vskip3mm
{\bf Proof.} The claims in $a)$ are clear, as well as the inequality and the monotonicity of $L_\ep^*$ 
with respect to $\ep$ in $b)$. Note that $L^*$ is a lower semi-continuous function, since it represents 
the supremum of a family of continuous (linear) functions. That is, for all $x \in \R^n$, 
$$
\liminf_{y \rightarrow x} L^*(y) = L^*(x).
$$
In particular, this property holds along every line. Since the functions $r \rightarrow L^*(rx)$ 
are non-decreasing and convex on the positive half-axis $r \geq 0$, we conlcude that 
$L^*((1-\ep)x) \uparrow L^*(x)$. 

Next, we have
\bee
L_\ep^*(x) 
 & = &
\sup_{y \in \R^n} \big[\left<x,y\right> - L(y) - \ep L_0(x)\big] \nonumber \\
 & = &
\sup_{y \in \R^n} \big[\left<(1-\ep)x,y\right> - L(y) + \ep (\left<x,y\right> - L_0(x))\big] \nonumber \\
 & \geq &
\sup_{y \in \R^n} \big[\left<(1-\ep)x,y\right> - L(y)\big] - 
\ep \sup_{y \in \R^n} \big[\left<x,y\right> - L_0(x)\big] \nonumber \\
 & = &
L^*((1-\ep)x) - \ep L_0^*(x).
\ene
Since $L_0^*(x)$ is finite, it follows that $\liminf_{\ep \rightarrow 0} L_\ep^*(x) \geq L^*(x)$ 
and hence $L_\ep^*(x) \uparrow L^*(x)$.

Since $L \leq L_\ep$, we have 
$$
\Big\{r>0: \int L_\ep(u/r)\,d\lambda \leq 1\Big\} \subset \Big\{r>0: \int L(u/r)\,d\lambda \leq 1\Big\},
$$
implying $\|u\|_{L_\ep} \geq \|u\|_L$. By the same argument, we get the monotonicity of $\|u\|_{L_\ep}$.

For the next claim, put $r = \|u\|_L$. We may assume that $r$ is finite and positive, so that
$\int L(u/r)\,d\lambda = 1$. Let $C = \int L_0(u/r)\,d\lambda$. Then
$\int L_\ep(u/r)\,d\lambda = 1 + C\ep$ and hence 
$$
\int L_\ep(u/r(1+C\ep))\,d\lambda \leq \frac{1}{1 + C\ep} \int L_\ep(u/r)\,d\lambda = 1,
$$
so that $\|u\|_{L_\ep} \leq r(1 + C\ep)$. It follows that 
$\limsup_{\ep \rightarrow 0} \|u\|_{L_\ep} \geq r$ and hence $\|u\|_{L_\ep} \downarrow r$.

For claim in $d)$, since $L_\ep^* \leq L^*$, we have 
$$
\Big\{r>0: \int L^*(u/r)\,d\lambda \leq 1\Big\} \subset \Big\{r>0: \int L^*_\ep(u/r)\,d\lambda \leq 1\Big\},
$$
implying $\|u\|_{L^*} \geq \|u\|_{L_\ep^*}$. By the same argument, we get the monotonicity of 
$\|u\|_{L_\ep^*}$.

For the last claim in $d)$, put $r = \|u\|_{L^*}$. There is nothing to prove, if $r = 0$. 
Assuming that $r$ is positive and finite, we have $\int L^*(u/r)\,d\lambda = 1$ and 
$\int L^*(u/r')\,d\lambda > 1$ whenever $0 < r' < r$. By $b)$ and applying the monotone 
convergence theorem, it follows that for such values of $r'$  
$$
\int L_\ep^*(u/r')\,d\lambda \, \uparrow \, \int L^*(u/r')\,d\lambda \quad {\rm as} \ \ \ep \downarrow 0.
$$
Since the last integral is greater than 1, this will be so for the first integral if $\ep$ is small enough,
and then $\|u\|_{L_\ep^*} \geq r'$, by the definition (8.1). Hence,
$$
\liminf_{\ep \rightarrow 0}\, \|u\|_{L_\ep^*} \geq r'.
$$
As $r' \in (0,r)$ was arbitrary, we may conclude that $\liminf_{\ep \rightarrow 0}\, \|u\|_{L_\ep^*} \geq r$ 
and thus $\|u\|_{L_\ep^*} \uparrow r$.

Finally, the remaining case $r = \infty$ means that $\int L^*(\alpha u)\,d\lambda > 1$ for all $\alpha > 0$.
Using $b)$ and applying the monotone convergence theorem, we have
$\int L_\ep^*(\alpha u)\,d\lambda > 1$ for all $\ep > 0$ small enough, in which case
$\alpha \|u\|_{L_\ep^*} > 1$. Hence,
$$
\liminf_{\ep \rightarrow 0}\, \|u\|_{L_\ep^*} > \frac{1}{\alpha}.
$$
As $\alpha > 0$ was arbitrary, we conclude that $\liminf_{\ep \rightarrow 0}\, \|u\|_{L_\ep^*} = \infty$ 
and thus $\|u\|_{L_\ep^*} \uparrow \infty$.
\qed

\vskip7mm
\section{{\bf Concavity of the Luxemburg Pseudo-Norm}}
\setcounter{equation}{0}

\vskip2mm
\noindent 
Let $L:\R^n \rightarrow [0,\infty]$ be a lower semi-continuous convex function such that 
$L(0)=0$, $L(x)>0$ for all $x \neq 0$, which is finite in some neighborhood of the origin.

One useful property of the Luxembourg norm is its quasi-concavity with respect to 
the measure $\lambda$. 
Given a measurable function $u:\Omega \rightarrow \R^n$, consider the functional 
\be
S(\lambda) = \|u\|_{L(\lambda)}
\en 
on the space of all probability measures $\lambda$ on $(\Omega,\frak M)$.

\vskip5mm
{\bf Proposition 11.1.} {\sl The functional $S$ is quasi-concave: Given probability measures 
$\lambda_1,\lambda_2$ on $\Omega$, for all $t_1,t_2 \geq 0$ such that $t_1 + t_2  = 1$, 
\be
S(t_1 \lambda_1 + t_2 \lambda_2) \geq \min\{S(\lambda_1),S(\lambda_2)\}.
\en
Moreover, if $L$ is positive homogeneous of order $p \geq 1$,
this functional is concave:
\be
S(t_1 \lambda_1 + t_2 \lambda_2) \geq t_1 S(\lambda_1) + t_2 S(\lambda_2).
\en
}

\vskip0mm
{\bf Proof.} Put $\lambda = t_1 \lambda_1 + t_2 \lambda_2$,
$r_1 = S(\lambda_1)$, $r_2 = S(\lambda_2)$.
One may assume that $0 < r_i < \infty$, so that $\int L(u/r_i)\,d\lambda_i = 1$.
Putting $r = \min(r_1,r_2)$, $\alpha_i = r_i/r$, and using $\alpha_i \geq 1$, we have
\bee
\int L(u/r)\,d\lambda
 & = &
t_1 \int L(\alpha_1 u/r_1)\,d\lambda_1 + t_2 \int L(\alpha_2 u/r_2)\,d\lambda_2 \\
 & \geq &
t_1 \int L(u/r_1)\,d\lambda_1 + t_2 \int L(u/r_2)\,d\lambda_2 \ = \ 1.
\ene
Hence, $\|u\|_{L(\lambda)} \geq r$, and (11.2) follows.

In the second claim, recall the definition (6.4) of the function $\Psi_L$. In the homogeneous
case, necessarily $\Psi_L(r) = r^p$ is convex on the half-axis $r \geq 0$.
Putting $r = t_1 r_1 + t_2 r_2$, by Jensen's inequality, we get
\bee
\int L(u/r)\,d\lambda
 & = &
t_1 \int L(\alpha_1 u/r_1)\,d\lambda_1 + t_2 \int L(\alpha_2 u/r_2)\,d\lambda_2 \\
 & \geq &
t_1 \int \Psi_L(\alpha_1) L(u/r_1)\,d\lambda_1 + t_2 \int \Psi_L(\alpha_2) L(u/r_2)\,d\lambda_2 \\
 & = &
t_1 \Psi_L(\alpha_1) + t_2 \Psi_L(\alpha_2) \\
 & \geq & 
\Psi_L(t_1 \alpha_1 + t_2 \alpha_2) \ = \ 1.
\ene
Hence, $\|u\|_{L(\lambda)} \geq r$, that is, (11.3).
\qed

\vskip5mm
The inequality (11.3) can be extended, as well as (11.2), to arbitrary (``continuous") convex mixtures 
of probability measures. Let us state such a relation for convolutions $\lambda * \kappa$ on 
the Euclidean space $\Omega = \R^m$, 
which we equip with the Borel $\sigma$-algebra. Note that any such convolution represents 
a convex mixture of shifts or translates of $\lambda$ with a mixing measure $\kappa$.

\vskip5mm
{\bf Proposition 11.2.} {\sl Suppose that $L$ is positive homogeneous of order $p \geq 1$.
Given a Borel measurable function $u:\R^m \rightarrow \R^n$,
for all probability measures $\lambda$ and $\kappa$ on $\R^m$,
\be
\|u\|_{L(\lambda * \kappa)} \, \geq \, \int \|u(x-y)\|_{L(\lambda(dx))}\,d\kappa(y).
\en
}

\vskip2mm
In the general case of a not necessarily homogeneous function $L$, the relations (11.3)-(11.4)
are extended under the $\Delta_2$-condition in a somewhat weaker form. Given 
$p_1 \geq p_0 \geq 1$,  define the quantity
\be
\gamma(p_1,p_0) = \sup \big\{a + b: a+br \leq \min(r^{p_1},r^{p_0}) \ {\rm for \ all} \ r \geq 0\big\}.
\en
In particular, $0 < \gamma(p_1,p_0) \leq 1$ and $\gamma(p_1,p_0) = 1$ if $p_1 = p_0$ (and only in this case).

According to Proposition 6.3, in presence of the two-sided $\Delta_2$-condition we have
\be
\Psi_L(r) = \inf_{x \neq 0} \frac{L(rx)}{L(x)} \geq \min(r^{p_1},r^{p_0})
\en
with
\bee
p_1 
 & = &
p_+ \, = \, \Phi_L'(1+) \, = \, \Psi_L'(1-), \\
p_0 
 & = &
p_- \, = \, \Phi_L'(1-) \, = \, \Psi_L'(1+).
\ene

\vskip3mm
{\bf Proposition 11.3.} {\sl Given probability measures 
$\lambda_1,\dots,\lambda_N$ on $\Omega$, for all $t_i \geq 0$ such that $t_1 + \dots + t_N  = 1$, we have
\be
S(t_1 \lambda_1 + \dots + t_N \lambda_N) \geq \gamma \sum_{i=1}^N t_i S(\lambda_i)
\en
with constant $\gamma = \gamma(p_+,p_-)$. As a consequence, if the function $u:\R^m \rightarrow \R^n$ is
Borel measurable, then for all probability measures $\lambda$ and $\kappa$ on $\R^m$,
$$
\|u\|_{L(\lambda * \kappa)} \, \geq \, \gamma \int \|u(x-y)\|_{L(\lambda(dx))}\,d\kappa(y).
$$
}

\vskip2mm
{\bf Proof.} Denote by $V(r)$ the maximal convex function majorized by 
$\min(r^{p_+},r^{p_-})$ on the positive half-axis $r \geq 0$. Then, by (11.5)-(11.6),
$$
\Psi_L(r) \geq V(r) \ {\rm for \ all} \ r \geq 0, \quad V(1) = \gamma.
$$

As before, assume that the pseudo-norms $r_i = S(\lambda_i)$ are positive and finite 
for all $i \leq n$, so that $\int L(u/r_i)\,d\lambda_i = 1$. Put $r = t_1 r_1 + \dots + t_N r_N$, 
$\alpha_i = r_i/r$. By Jensen's inequality, 
\bee
\int L(u/r)\,d\lambda
 & = &
t_1 \int L(\alpha_1 u/r_1)\,d\lambda_1 + \dots + t_N \int L(\alpha_N u/r_N)\,d\lambda_N \\
 & \geq &
t_1 \int \Psi_L(\alpha_1) L(u/r_1)\,d\lambda_1 + \dots + t_N \int \Psi_L(\alpha_N) L(u/r_N)\,d\lambda_N \\
 & = &
t_1 \Psi_L(\alpha_1) + \dots + t_N \Psi_L(\alpha_N) \\
 & \geq &
t_1 V(\alpha_1) + \dots + t_N V(\alpha_N) \\ 
 & \geq & 
V(t_1 \alpha_1 + \dots + t_N \alpha_N) \, = \, V(1) \, = \, \gamma.
\ene
Since $0 < \gamma \leq 1$, it follows that
$$
\int L\Big(\frac{u}{\gamma r}\Big)\,d\lambda \geq \frac{1}{\gamma}\int L\Big(\frac{u}{r}\Big)\,d\lambda \geq 1.
$$
Hence, $\|u\|_{L(\lambda)} \geq \gamma r$, that is, (11.7).
\qed

\vskip5mm
{\bf Acknowledgement.} We would like to thank Lech Maligranda for pointing out
missing references related to the history of Luxemburg and Orlicz norms.

\end{document}